\documentclass[final]{siamltex}
\usepackage{amsfonts,amsmath,enumerate,hyperref,bm}

\newcommand{\eps}{\varepsilon}
\newcommand{\cU}{{\cal U}}

\newcommand{\cE}{{\cal E}}

\newcommand{\Tmix}{T_{\mathrm{mix}}}
\newcommand{\reals}{\mathbb{R}}

\newcommand{\Prob}{\ensuremath{{\bf Pr}}}

\newcommand{\E}{\ensuremath{{\bf E}}}
\newcommand{\Ent}{\ensuremath{{\bf Ent}}}
\newcommand{\LS}{\alpha_{\mathrm{E}}}

\newtheorem{remark}[theorem]{Remark}
\newtheorem{claim}[theorem]{Claim}

\title{Reconstruction for Colorings on Trees}
\author{Nayantara Bhatnagar\thanks{Department of Statistics,
    University of California, Berkeley, {\tt
      nayan@stat.berkeley.edu}. Research supported in part by
      NSF grant CCF-0455666 and DOD ONR grant
    N0014-07-1-05-06.}
\and Juan Vera\thanks{Department of Management Sciences,
University of Waterloo, {\tt jvera@uwaterloo.ca}}
\and
Eric Vigoda\thanks{College of Computing,
Georgia Institute of Technology, {\tt vigoda@cc.gatech.edu}.
Research supported in part by NSF grant CCF-0455666.}
\and
Dror Weitz \thanks{
Tel Aviv, Israel,
{\tt dror@drorweitz.com}.}
}

\begin{document}

\title{Reconstruction for Colorings on Trees}

\maketitle

\begin{abstract}
  Consider $k$-colorings of the complete tree of depth $\ell$ and branching
  factor $\Delta$.  If we fix the coloring of the leaves, for what range of $k$
  is the root uniformly distributed over all $k$ colors (in the limit
  $\ell\rightarrow\infty$)?  This corresponds to the threshold for
  uniqueness of
  the infinite-volume Gibbs measure.  It is straightforward to show the
  existence of colorings of the leaves which ``freeze'' the entire tree when
  $k\le\Delta+1$. For $k\geq\Delta+2$, Jonasson proved the root is ``unbiased''
  for any fixed coloring of the leaves and thus the Gibbs measure is unique.
  What happens for a {\em typical} coloring of the leaves?  When the
  leaves have
  a non-vanishing influence on the root in expectation, over random
  colorings of the leaves, reconstruction is said to hold.
  Non-reconstruction is equivalent
  to extremality of the free-boundary Gibbs measure.  When
  $k<\Delta/\ln{\Delta}$, it is straightforward to show that reconstruction is
  possible (and hence the measure is not extremal).

  We prove that for $C>1$ and $k =C\Delta/\ln{\Delta}$, non-reconstruction
  holds, i.e., the Gibbs measure is extremal. We prove a strong form of
  extremality: with high probability over the colorings of the leaves the
  influence at the root decays exponentially fast with the depth of the tree.
  Closely related results were also proven recently by Sly. The above
  strong form of
  extremality implies that a local Markov chain that updates constant sized
  blocks has inverse linear entropy constant and hence $O(N\log N)$ mixing time
  where $N$ is the number of vertices of the tree.  Extremality on trees and
  random graphs has received considerable attention recently since it may have
  connections to the efficiency of local algorithms.
\end{abstract}

\begin{keywords} 
Reconstruction, Random colorings, Extremality of Gibbs measure
\end{keywords}

\begin{AMS}
60J05, 60K35
\end{AMS}

\pagestyle{myheadings}
\thispagestyle{plain}

\section{Introduction}

Reconstruction on trees arises naturally in a variety of contexts.
In evolutionary biology it is intimately related with the efficiency of
inferring phylogenetic ancestors \cite{Mossel:phylo,CMR}.
It appears in communication theory in the study of noisy computation
\cite{EKPS}.
Most closely related to the origins of this work, for spin systems
from statistical physics, the threshold for reconstruction is equivalent to
the threshold for extremality of the infinite-volume
Gibbs measure induced by free-boundary conditions
(see, e.g., \cite{Geo}), and is closely connected
to mixing properties of the single-site update Markov chain known
as the Glauber dynamics \cite{BKMP,MSW2}.
In fact, Berger et al \cite{BKMP} showed that $O(N\log{N})$ mixing
time of the Glauber
dynamics implies extremality.
Recent work of Krzakala et al \cite{KMRSZ} used statistical physics
methods to study extremality for random colorings and random $k$-SAT,
and suggests there are interesting connections between
thresholds for extremality and efficiency of local algorithms.

Our general setup is a complete tree with branching factor
$\Delta$ and depth $\ell$, which
we will denote as $T_\ell=(V_\ell,E_\ell)$.
(To clarify, the maximum degree of this tree is $\Delta+1$.)
A $k$-coloring is an assignment
$\sigma:V_\ell\rightarrow [k]=\{1,\dots,k\}$ where
$(v,w)\in E_\ell $ implies $\sigma(v)\neq\sigma(w)$.
Let $\mu_\ell$ denote the uniform measure over $k$-colorings of
the tree.  Finally, we denote the leaves of the tree $T_\ell$ by
$L_\ell$. When the height of the tree is clear from the context, we
will drop the subscript and just use $L$.

In the reconstruction problem, we are asking whether a typical coloring of
the leaves influences the conditional measure at the root.  In words,
we choose a random coloring
of the tree, fix the induced coloring for the leaves, and then choose a random
coloring of the internal vertices consistent with the leaves.
If the root has some non-vanishing bias to any color then we say
reconstruction is possible. Here is a precise definition.

\begin{definition}
The reconstruction problem for $T_\ell$ is {\em solvable} if there
exists $c \in [k]$ such that
\[
\lim_{\ell \rightarrow \infty}
\E_{\tau\sim\mu_\ell}
\left[ \Big|    \mu_\ell\big(\sigma(\mathrm{root})=c \, \vert \,
  \sigma(L)=\tau(L)\big)  -1/k \Big| \right]
    > 0.
\]
\end{definition}
Non-reconstruction is  equivalent to extremality of
the corresponding Gibbs measure of the infinite
tree (i.e., it can not be expressed as a convex combination of other
Gibbs measures), see, e.g., \cite{Geo}.

In the corresponding reconstruction problem for the Ising model, the
reconstruction threshold has been
precisely established \cite{BRZ}.  Recent works have improved
bounds on the reconstruction threshold for the case of the hard-core
model on weighted independent sets \cite{BWextrem,Martin,MSW2}.
More recently, Borgs et al \cite{BCMR} showed that for binary
asymmetric channels, the Kesten-Stigum eigenvalue
bound is tight for reconstruction non-solvability in the case of limited
asymmetry.

The problem of determining the reconstruction threshold for
colorings on the tree was discussed in \cite{BW,MM,Mossel}.
The threshold for uniqueness of the Gibbs measure on trees is now well-known.
In contrast to extremality, uniqueness of the Gibbs measure requires that for
every fixed coloring of the leaves, the root is not biased to any
color.  Jonasson \cite{Jonasson}
proved that for $k\geq\Delta+2$ the Gibbs measure is unique \cite{Jonasson},
whereas for $k=\Delta+1$ there are ``frozen'' boundary conditions in the
sense that the root has only one possible color, see \cite{BW} for more on the
topic of uniqueness for colorings.

It was known \cite{MP} that if $k < (1-\varepsilon)\Delta/\ln{\Delta}$ for
any $\varepsilon>0$ then reconstruction is solvable, and more
precisely as established in \cite{Semerjian}, reconstruction is
solvable for $\Delta \ge k(\ln k + \ln\ln k + 1 +o(1))$.
Recent works \cite{KMRSZ,ZK} used the
so-called replica-symmetry breaking method to conjecture that the
reconstruction
threshold for $k$-colorings of the tree $T_\ell$ is at
$\Delta =
k(\ln k + \ln\ln k + 1 +o(1))$.   From an algorithmic
perspective it is interesting that these bounds on the degree are
close to the current best
results for efficiently coloring Erd\"{o}s-R\'{e}nyi random graphs
$G(n,\Delta/n)$ \cite{AM}.

In this work we prove that the reconstruction threshold is at
$\frac{\Delta}{\ln{\Delta}}(1+o(1))$.
In particular, we show that for $C>1$ and $k>
C\Delta/\ln\Delta$ reconstruction is not possible, with the bias of an average
boundary condition decaying exponentially fast with the height of the tree.
\begin{theorem}
\label{th:non-reconstruction}
  Let $C>1$ and $\eps = \min\{C-1,1/3\}$. There are constants $\alpha(\eps)>0$,
  $\Delta_0 = \Delta_0(\eps)$, and $\ell_0 = \ell_0(\eps)$ such that for all
  $\Delta > \Delta_0$, every $k > C\Delta/\ln\Delta$, and every color $c \in
  [k]$, for $\ell>\ell_0$,
\begin{eqnarray*}
\E_{\tau \sim \mu_\ell}\left[ \Big|\mu_\ell\big(\sigma(\mathrm{root})=c
  \,\vert\, \sigma(L)=\tau(L)\big) -1/k \Big| \right]
\leq
    \Delta^{-\alpha\ell}.
\end{eqnarray*}
\end{theorem}
We also prove a high concentration result, where a much stronger tail
bound (on the fraction of biasing boundary conditions) is established.
\begin{theorem}
\label{th:concentration}
Let $C>1$ and $\eps = \min\{C-1,1/3\}$. There are constants
$\alpha(\eps),\alpha'(\eps)>0$, $\Delta_0 = \Delta_0(\eps)$ and
$\ell_0 = \ell_0(\eps)$
such that for all $\Delta > \Delta_0$, every $k > C\Delta/\ln\Delta$,
and every color $c \in [k]$, for $\ell>\ell_0$,
\begin{eqnarray*}
\Prob_{\tau\sim\mu_\ell}\left[\Big| \mu_\ell\big(\sigma(\mathrm{root})
  = c
  \,\vert\,
\sigma(L) = \tau(L) \big) - 1/k\Big|
> \Delta^{-\alpha\ell}\right]
 \leq e^{-\Delta^{\alpha'\ell}}.
\end{eqnarray*}
\end{theorem}

\begin{remark}[A Note on Comparison with Previous Work]
An earlier version of this paper proved the above results for $C>2$,
which was the first work giving the right order for the threshold for
non-reconstruction \cite{BVV}. Stated in terms to make comparison with
other work easier, it was shown there that non-reconstruction holds when
$\Delta \leq \frac{1}{2} k \ln k - o(k \ln k)$. The bound was
sharpened by Sly in \cite{Sly}, using independent methods, and showing that
non-reconstruction holds when $\Delta \leq k[\ln{k} + \ln\ln{k} + 1 -\ln{2} -
o(1)]$. Subsequent to this the authors of this paper refined the approach of
\cite{BVV} to obtain the current results, showing non-reconstruction
when $\Delta \leq  k \ln k - o(k \ln k)$. Though the bound on the
reconstruction threshold in \cite{Sly} is more precise, we feel our
results are still of independent interest for the following reasons.

The techniques used in the current work (an extension of
those appearing in \cite{BVV}, utilizing ideas from \cite{MSW1})
are independent and different from the
analytical approach of \cite{Sly}. There may be some insight to be
gained from this proof approach. Secondly, as a consequence
of the strong concentration result stated above, using results of
Martinelli et al \cite{MSW1,MSW2}, it follows that a local Markov
chain (in each step the colors of a constant number of vertices are updated)
has inverse linear entropy constant and hence $O(N\log{N})$
mixing time.
\end{remark}

We state the relevant definitions before formally stating the
result for the dynamics. We consider dynamics on the tree $T_n$, and
analyze its performance as a
function of the volume of the tree $N=|T_n|$. In each step the dynamics updates
the colors of a block of vertices. Specifically, for a fixed parameter $\ell$,
for every $v\in T_n$ let $B_{v,\ell}$ be the subtree (of $T_n$) of depth $\ell$
rooted at $v$, i.e., $B_{v,\ell}$ includes all vertices in $T_n$ that are
descendants of $v$ (including $v$) and whose distance from $v$ is at most
$\ell$.  We analyze the dynamics $\mathcal{M}_\ell$, which makes
heat-bath updates
in a random block $B_{v,\ell}$. Specifically, let the state space $\Omega$ be
the set of proper $k$-colorings of $T_n$. From $X_t\in\Omega$, the transition
$X_t\rightarrow X_{t+1}$ is defined as follows:
\begin{itemize}
\item
Choose a vertex $v$ uniformly at random from $T_n$.
\item
For all $w\not\in B_{v,\ell}$, set $X_{t+1}(w)=X_t(w)$.
\item
Choose $X_{t+1}(B_{v,\ell})$ uniformly at random from the
set of $k$-colorings consistent with $X_{t+1}(T_n\setminus B_{v,\ell})$.
\end{itemize}

The case when $\ell=0$ is the well-studied
single-site dynamics known as the Glauber dynamics.
It is easy to see that for $k\geq 4$
and for all $\ell$, $\mathcal{M}_\ell$ is an ergodic Markov chain with unique stationary
distribution $\pi$ uniformly distributed over $\Omega$.
The central question is the mixing time defined as:
\[ \Tmix = \max_{X_0} \min\{t: \| P^t(X_0,\cdot) - \pi\| \leq 1/2e \} \]
where $P^t(X_0,\cdot)$ is the distribution after $t$ steps of the dynamics
starting from coloring $X_0$ and $\|\cdot\|$ denotes variation distance.
The choice of the constant $1/2e$ implies
that variation distance $\leq\eps$ can be achieved after $\leq\lceil\ln{1/\eps}\rceil\Tmix$
steps \cite{Aldous}.

We bound the mixing time via the entropy constant.
Let $f:\Omega\rightarrow\reals$ be an arbitrary test function.
The entropy of~$f$ is
\[
\Ent(f) := \Ent_\pi(f) = \E_\pi[f\ln(f)] - \mathbf{E}_\pi(f)\ln[\mathbf{E}_\pi(f)],
\]
and the entropy constant of the Markov chain is defined as
\[
\LS = \inf_{f\geq 0} \frac{\Ent(f) - \Ent(Pf)}{\Ent(f)}
\]
where the infimum is over non-constant functions $f$.  The entropy constant,
which bounds the rate of decay of entropy, provides a good bound on the mixing
time.  In particular, standard results (see,
e.g., \cite{DSc,FK,MT}) imply:
\[
\Tmix \leq O[\LS^{-1}\log\log(\pi_{\min}^{-1})] \le O(\LS^{-1}\log{N}).
\]
We prove the following result about the dynamics.
\begin{theorem}
\label{th:dynamics}
Let $C>1$ and $\eps = \min\{C-1,1/3\}$. There are constants $\Delta_0 =
\Delta_0(\eps)$, and $\ell = \ell(\eps)$ such that for all $\Delta > \Delta_0$,
every $k > C\Delta/\ln\Delta$, and all $n$, the entropy constant of the dynamics
$\mathcal{M}_\ell$ on $T_n$ satisfies
\[ \LS \geq \frac{1}{4N}
\]
and consequently the mixing time satisfies
\[ \Tmix = O(N\log{N}).
\]
\end{theorem}

It is an interesting open problem to prove $O(N\log{N})$ mixing time of the
Glauber dynamics for the same range of colors as the above theorem.
The best known results for the mixing time of the Glauber dynamics on the
tree are
Martinelli et al \cite{MSW2} who proved $O(N\log{N})$ mixing time
when $k\geq\Delta+3$ for arbitrary boundary conditions, and
Hayes et al \cite{HVV} proved polynomial in $N$
mixing time (specifically, $O^*(N^4)$ mixing time)
when $k>100\Delta/\log{\Delta}$ for any planar graph.
For complete trees, Goldberg et al \cite{GJK} and Lucier et al \cite{LMP} recently showed
the mixing time is polynomial for any fixed $k$ and $\Delta$.

Finally, a related problem is the reconstruction
threshold on Erd\"{o}s-R\'{e}nyi random graphs $G(n,p)$ with $p=\Delta/n$,
so-called sparse random graphs. Recent work of Gerschenfeld and
Montanari \cite{GM} gives sufficient conditions under which
extremality on the tree is equivalent to extremality for sparse random
graphs. They showed that
the conditions are satisfied for the $q$-state
Potts model at all temperatures except zero temperature, which is the case of
proper colorings. Subsequently in \cite{MRT}, the conditions were
also verified in the case of proper
colorings.

\section{Proof outline and outline of paper}
\label{sec:outline}
Our proof is divided into two major parts.
In the first part, which is contained in Section \ref{sec:leafcolorings},
we prove that with very high probability over the colorings of the
leaves, the root is not too highly biased in favor of any color.  Roughly,
with probability $\geq 1-e^{-\Delta^{\Omega(\ell)}}$ over random
colorings $X$ of the leaves, for every
color $c$,
\[ \mu_\ell\left(\sigma(\text{root})=c \,\vert\, \sigma(L)=X\right)
\leq \Delta^{-\eps/2}
\]
where $\eps$ is as defined in Theorem~\ref{th:concentration}.  We prove this
statement by analyzing the same recurrences as used by Jonasson \cite{Jonasson}
in his proof for uniqueness on the tree when $k\geq\Delta+2$.  His recurrences
express the marginal distribution at the root of a tree of size $\ell$ in terms
of the marginals for trees of size $\ell-1$.  The difficulty in our setting is
that when $k<\Delta+2$ it is unclear if the recurrences converge to the uniform
distribution as a fixed point.

In the second part of the proof (which appears in Section \ref{sec:convergence})
we use a two stage coupling similar to that used
in \cite{MSW1} together with the bound on the maximum probability of a
color from the first part to establish the decay of correlation stated in
Theorem~\ref{th:non-reconstruction}. We then add arguments (taken
from \cite{MSW1}) to get the stronger Theorem~\ref{th:concentration}.

Finally, in Section \ref{sec:dynamics} we prove Theorem \ref{th:dynamics}
establishing fast convergence of the block dynamics for constant sized blocks.

\section{Unbiasing leaf colorings}
\label{sec:leafcolorings}
In this section we show that for most colorings of the leaves
$X$, at all vertices far enough from the leaves, the color
conditioned on~$X$ has ``sufficient'' randomness.

We call a vector $X \in \{\star,1,\dots,k\}^{\Delta^h}$ a
{\it (partial) coloring of the leaves} $L_h$. The notation $X_i =
\star$ is used to denote that
the $i$-th leaf is not assigned a color.
We say that $X$ is {\em allowed} if there exists a coloring $\sigma$
of the tree $T_h$
which is consistent with $X$, i.e., for every leaf $z$, if
$X(z)\in\{1,\dots,k\}$ then
$\sigma(z)=X(z)$.

Given an allowed coloring $X$ of $L_h$  and $c \in [k]$ we define
  $P_h(X,c)$ inductively by
\[
P_0(X,c) = \begin{cases}
1 &\text{if }X = c\\
1/k &\text{if }X = \star\\
0 &\text{oherwise}
\end{cases}
\hspace{.7in}
P_h(X,c)  =  \frac{\displaystyle\prod_{i=1}^\Delta
  (1-P_{h-1}(X_i,c))}
{\displaystyle\sum_{d \in [k]} \displaystyle\prod_{i=1}^\Delta
  (1-P_{h-1}(X_i,d))}
\]
where $X = (X_1,\dots,X_\Delta)$ with each $X_i \in
\{\star,1,\dots,k\}^{\Delta^{h-1}}$. An inductive argument shows that
$P_h(X,c)$ is well-defined for all allowed colorings $X$ and all colors
$c$.

A simple computation shows that for allowed colorings $X$, $P_h(X,c)$ is,
in fact, the probability the root is colored $c$ conditioned on $X$ at
the leaves. Note that conditioning on a partial coloring $X$ simply
means conditioning on the vertices which are assigned colors in $X$.
\begin{lemma}
\label{lem:mu=p}
For all $h$,  and all allowed colorings $X$ of $L_h$, and all $c$,
\[
\mu_{h}(\sigma(root) = c \, \vert \, \sigma(L_h) = X) = P_h(X,c).
\]
\end{lemma}
\begin{proof}
The proof is by counting the appropriate sets of colorings and induction on
$h$. Let $\Omega_h(X,c)$ be the number of colorings of a tree of height $h$
where the root is colored $c$ and the coloring is consistent with~$X$. Then,
\begin{eqnarray*}
\mu_{h}(\sigma(root) = c|\sigma(L_h) = X) & = &
  \frac{\Omega_h(X,c)}{\displaystyle\sum_{d \in k}
  \Omega_h(X,d)}
=
  \frac{\displaystyle\prod_{i=1}^\Delta \sum_{f \neq c}
  \Omega_{h-1}(X_i,f)}{\displaystyle\sum_{d
  \in k} \displaystyle\prod_{i=1}^\Delta \sum_{f \neq d}
  \Omega_{h-1}(X_i,f)   }.
\end{eqnarray*}
By the assumption that $X$ extends to a coloring of $T_h$, it must
be that each $X_i$ extends to a coloring of the tree $T_{h-1}$
hence, $\prod_i \sum_{f} \Omega_{h-1}(X_i,f) \neq 0$. Dividing
the numerator and denominator by this factor, we obtain
\begin{eqnarray*}
\frac{\displaystyle\prod_{i=1}^\Delta
  1-\frac{\Omega_{h-1}(X_i,c)}{\sum_{f} \Omega_{h-1}(X_i,f)
  }
}
{\displaystyle\sum_{d
  \in k} \displaystyle\prod_{i=1}^\Delta
  1-\frac{\Omega_{h-1}(X_i,d)}{\sum_{f}
  \Omega_{h-1}(X_i,f)}}
=
\frac{\displaystyle\prod_{i=1}^\Delta
  1-P_{h-1}(X_i,c)
}
{\displaystyle\sum_{d
  \in k} \displaystyle\prod_{i=1}^\Delta
  1-P_{h-1}(X_i,d)
} =P_h(X,c).
\end{eqnarray*}
\end{proof}

Henceforth, we assume that $C>1$ and $k = C\Delta/\ln{\Delta}$.  Define the
parameter $\eps = \eps(C) = \min\{C-1,1/3\}$. Notice that $k >
(1+\eps)\Delta/\ln{\Delta}$.  We also assume $\Delta_0(\eps)$ is a large
constant depending only on~$\eps$.

We now give a recursive definition of a coloring of the leaves being {\em
  unbiasing}. In the base case, let $X\in\{\star,1,\ldots,k\}^\Delta$ be a
(partial) coloring of $L_1$. We say that~$X$ is {\em unbiasing} if and only if
at least $\Delta^{\eps/2}$ colors do not appear in $X$. For $\ell>1$ let
$X=(X_1,\ldots,X_\Delta)$ be a coloring of $L_\ell$ where $X_i$ is the coloring
of the leaves of the subtree rooted at the $i$-th child of the root.  We say
that $X$ is unbiasing if and only if at most $\Delta^{1-\eps}$ of the $X_i$ are
not unbiasing.

Indeed, given an unbiasing coloring of the leaves, the color at root cannot be
too biased.
\begin{lemma}
\label{lem:unbiasing}
  For any $\ell\geq 1$, any unbiasing coloring~$X$ of $L_\ell$, and any
  color~$c$,
  $$P_\ell(X,c) \leq \Delta^{-\eps/2}.$$
\end{lemma}

The main result in this section is that if~$X$ is a random coloring of the
leaves~$L_\ell$, then~$X$ is unbiasing with very high probability.  We use
$X\sim\mu_\ell(\sigma(L))$ to denote that $X$ is a random coloring of the leaves
of the tree $T_\ell$, i.e., choose $\sigma$ from $\mu_\ell$ and let
$X=\sigma(L)$.
\begin{theorem}
\label{th:unbiasing-tail}
  For all $\ell\geq 1$,
  \[
  \Prob_{X\sim\mu_{\ell}(\sigma(L))}\left[ X \ \text{is unbiasing}\right] \geq
  1-e^{-\Delta^{\frac{\ell-(1-\eps)}{2}}}.
  \]
\end{theorem}

For a given~$\ell$, let $h(v)$ be the distance (or height) of the
vertex~$v$ from the leaves $L_\ell$, and $X_v$ be the restriction of the
coloring $X$ to those leaves that are in the subtree rooted at~$v$.
We call $X$
{\em highly unbiasing} if for every $v$ with $h(v)\ge \eps\ell$, $X_v$ is
unbiasing.
\begin{corollary}
\label{cor:highly-unbiasing-tail}
  For all $\ell > \ell_0(\epsilon)$,
  \[
  \Prob_{X\sim\mu_{\ell}(\sigma(L))}\left[X \ \text{is highly unbiasing}\right] \geq
  1-e^{-\Delta^{\varepsilon \ell/3}}.
  \]
\end{corollary}
\begin{proof}
  This follows from Theorem~\ref{th:unbiasing-tail} using a simple union bound
  once we notice that there are $O(\Delta^{(1-\eps)\ell})$ vertices~$v$ with
  $h(v)\ge\eps\ell$.
\end{proof}

The rest of this section is dedicated to proving Lemma~\ref{lem:unbiasing} and
Theorem~\ref{th:unbiasing-tail}.

\subsection{Properties of partial colorings of the leaves}
We now show that an unbiasing coloring of the leaves indeed implies a
 not-too-biased color at the root.
\begin{proof}[Proof of Lemma \ref{lem:unbiasing}]
  The proof is by induction on~$\ell$. Let $X$ be an unbiasing coloring of the
  leaves $L_\ell$. We need to show that $P_\ell(X,c) \le \Delta^{-\eps/2}$ for
  every color $c\in[k]$.  For the base case when $\ell=1$, by definition, there
  are at least $\Delta^{\eps/2}$ colors not appearing in~$X$, and hence
  $P_1(X,c) \le \Delta^{-\eps/2}$ for every color~$c$. For $\ell>1$, write
  $X=(X_1,\ldots,X_\Delta)$, where $X_i$ is the coloring of the leaves of the
  subtree rooted at the $i$-th child of the root.  By definition, at most
  $\Delta^{1-\eps}$ of the $X_i$ are not unbiasing. By induction, if $X_i$ is
  unbiasing then for all colors~$c$, $P_{\ell-1}(X_i,c) \le \Delta^{-\eps/2}$.
  If, however, $X_i$ is not unbiasing then there are at most $\Delta^{\eps/2}$
  colors $c$ for which $P_{\ell-1}(X_i,c) \ge \Delta^{-\eps/2}$ because by
  Lemma~\ref{lem:mu=p} we have $\sum_{d \in [k]} P_{h-1}(X_i,d) = 1$.
  Therefore, there is a set $G \subseteq [k]$ of size at least $k -
  \Delta^{1-\eps}\Delta^{\eps/2} = k-\Delta^{1-\eps/2}$ such that for all $d \in
  G$, and for every $i \in [\Delta]$, $P_{\ell-1}(X_i,d) \leq \Delta^{-\eps/2}$.
  Hence, we have
    \[
    P_\ell(X,c)  =  \frac{\displaystyle\prod_{i=1}^\Delta (1-P_{\ell-1}(X_i,c))}
    {\displaystyle\sum_{d \in [k]} \displaystyle\prod_{i=1}^\Delta (1-P_{\ell-1}(X_i,d))}
    \leq  \frac{1} {\displaystyle\sum_{d \in G} \displaystyle\prod_{i=1}^\Delta
      (1-P_{\ell-1}(X_i,d))}.
     \]
    Using the arithmetic-geometric mean inequality,
  \[
   P_\ell(X,c)  \leq  \frac{1} {|G| \displaystyle\prod_{d \in G} \displaystyle\prod_{i=1}^\Delta
      (1-P_{\ell-1}(X_i,d))^{1/|G|}}
   \leq  \frac{1}{|G|}\exp\left(\frac{1+\Delta^{-\eps/2}}{|G|}\displaystyle\sum_{d \in G}
        \displaystyle\sum_{i=1}^\Delta P_{\ell-1}(X_i,d) \right)
   \]
   where the second inequality uses the fact that if $x<\delta <
  1/10$ then $1-x\geq \exp(-(1+\delta)x)$, and that $\Delta>\Delta_0(\eps)$. Using induction and
  plugging in the lower bound $|G| \geq k - \Delta^{1-\eps/2} \geq  (1-\Delta^{-\eps/3})k $ (for
  $\Delta$ large enough), we obtain
   \[
     P_\ell(X,c) \leq  \frac{1}{|G|}\exp\left(\frac{(1+\Delta^{-\eps/2})\Delta}{|G|}\right)
   \leq
    \frac{\Delta^{\frac{1 + \Delta^{-\eps/2}}{(1+\eps)(1-\Delta^{-\eps/3})}}}
    {k(1-\Delta^{-\eps/3})}\leq
    \Delta^{-\eps/2}.
    \qquad
  \]
\end{proof}

Before we go on to bound the probability that a random coloring of the leaves is
unbiasing (i.e., prove Theorem~\ref{th:unbiasing-tail}), we state some basic
lemmas which give some independence, allowing a recursive solution.
\begin{lemma}
\label{lem:identical}
  Let $X\sim\mu_\ell(\sigma(L))$.  For any $0\leq h \leq \ell$, let
  $w_1,\dots,w_{\Delta^{\ell-h}}$ denote the vertices which are at depth
  $\ell-h$ from the root.  Let $X=\left(X_{w_1},\dots,X_{w_{\Delta^{\ell-h}}}\right)$ where
  $X_{w_i}$ is the coloring of the leaves of the subtree rooted at $w_i$.  Then
  the $X_{w_i}$ are identically distributed as $\mu_h(\sigma(L_h))$.
\end{lemma}
\begin{proof}
  Consider the following recursive method for constructing a random coloring of
  the leaves: choose the color of the root uniformly at random from $k$ colors
  and independently choose the color of each child from the $k-1$ remaining
  colors. Now consider a vertex $v$ at height $h$ from the leaves.  Since each
  color is equally likely to appear at $v$, the distribution of the coloring at
  the leaves of the subtree rooted at $v$ is identical to the distribution over
  colorings of leaves in a random coloring of $T_{h}$.
\end{proof}

\begin{lemma}
\label{lem:independenceofgood}
  Let $X = (X_1, \cdots, X_\Delta)\sim\mu_\ell(\sigma(L))$ where $X_i\in
  [k]^{\Delta^{\ell-1}}$ is $X$ restricted to the subtree rooted at the $i$-th
  child of the root.  Let $\cU_i$ denote the event that $X_i$ is unbiasing.  The
  events $\cU_i, 1\le i\le\Delta$ are independent.
\end{lemma}
\begin{proof}
  The proof follows from the Markovian property of the Gibbs measure (the
  configurations on the subtrees are independent of each other once we condition
  on a spin at the root), and the fact that the event of being unbiasing is
  symmetric with respect to colors. Formally, for $I \subseteq [\Delta]$, let $\cU_I =
  \bigcap_{i \in I} \cU_i$.  It is enough to show that if $j \notin
  I$, then $\mu_\ell(\cU_j | \cU_I) = \mu_\ell(\cU_j)$. Now,
  \begin{eqnarray*}
  \mu_\ell(\cU_j | \cU_I) &  = & \sum_{c\in[k]} \mu_\ell(\sigma(root) = c \,|\, \cU_I)
  \cdot \mu_\ell(\cU_j \,|\, \sigma(root) = c \,,\, \cU_I)\\
  & = & \sum_{c\in[k]} \mu_\ell(\sigma(root) = c \,|\, \cU_I)\cdot
  \mu_\ell(\cU_j \,|\, \sigma(root) = c)\\
  & = & \sum_{c\in[k]} \mu_\ell(\sigma(root) = c) \cdot
  \mu_\ell(\cU_j \,|\, \sigma(root) = c)\\
  & = & \mu_\ell(\cU_j),
  \end{eqnarray*}
  where the second equality follows from the Markovian property of the Gibbs
  measure and the third equality uses the fact that the event $\cU_I$ is
  symmetric with respect to the color at the root.
\end{proof}

We are now ready to complete the proof of the main theorem in this section.
\begin{proof}[Proof of Theorem \ref{th:unbiasing-tail}]
  Let $q_\ell = \Prob_{X\sim\mu_\ell(\sigma(L))}\left[X \mbox{ is not unbiasing}\right].$
  We need to show that
  $$
  q_\ell \leq e^{-\Delta^{(\ell-(1-\eps))/2}},
  $$
  which we will do by induction on~$\ell$. For $\ell>1$ let the random
  coloring $X = (X_1,\dots,X_\Delta)$, where $X_i$ is the random coloring of the
  leaves of the subtree rooted at the $i$-th child of the root. Assuming
  $\Delta\geq\Delta_0(\eps)$ (so that $e^{-\Delta^{\eps/2}} < \Delta^{-2}$) we
  have
  \begin{align*}
    q_\ell & =  \Prob_{X\sim\mu_\ell(\sigma(L))}\left[> \Delta^{1-\eps}
    \text{ of the } X_i \text{ are not unbiasing}\right]
     &\text{(by definition)}
     \\
    & \leq \Delta^{\Delta^{1-\eps}} \cdot (q_{\ell-1})^{\Delta^{1-\eps}}
&   \mbox{(by Lemmas \ref{lem:independenceofgood} and \ref{lem:identical})}  \\
    & \leq (\Delta e^{-\Delta^{(\ell-(2-\eps))/2}})^{\Delta^{1-\eps}}
&        \text{(by induction)} \\
    & \leq e^{-\Delta^{(\ell-(1-\eps))/2}}.
  \end{align*}
  In the base case of $\ell=1$ the colors at the leaves are chosen independently
  (from $k-1$ colors, say from $\{1,\cdots,k-1\})$ once the color at the
  root is fixed. We have to upper
  bound the probability that there are less than $\Delta^{\eps/2}$ unused
  colors. Let $S$ denote the set of unused colors. Let $x_{i,c}$ be a
  0-1 valued random variable which is $1$
  if the color $c$ is chosen at the vertex $i$. Then the number of
  unused colors is given by
\begin{align*}
|S|  = \sum_{c=1}^{k-1} \prod_{i=1}^\Delta (1-x_{i,c})
\end{align*}
By linearity of expectation and independence over the children, we have
\begin{align*}
\E[|S|]  = \sum_{c=1}^{k-1} \prod_{i=1}^\Delta (1-\E[x_{i,c}]) =
k\left(1-\frac{1}{k-1}\right)^\Delta \geq 10 \Delta^{\varepsilon/2},
\end{align*}
where the last inequality follows by using the fact that for
  $\Delta>\Delta_0(\eps)$, $\Delta \le
  (1-\frac{2\eps}{3})(k-1)\ln (k-1)$ and, for $x<\delta$, $1-x
  \geq \exp(-(1+\delta)x)$.
It is not difficult to verify (for example by induction) that for $Y_c
  := \prod_i(1-x_{i,c})$, the variables $\{e^{t Y_c}\}_{c \in [q]}$
  have negative
  covariance for any $t>0$. Hence, we can apply
  the usual Chernoff bounds, (see, e.g.,
  the proof of Proposition 7 of
  \cite{DR}), to show that $\Prob[|S| < \Delta^{\varepsilon/2}] = q_1 <
  e^{-\Delta^{\eps/2}}$.
\end{proof}

\section{Convergence}
\label{sec:convergence}
In this section we complete the proof of Theorem~\ref{th:non-reconstruction},
i.e., establish exponential decay of correlations for the measure. We do that
using coupling arguments where disagreements percolate down and up the tree,
along the lines suggested in \cite{MSW1}. The main addition in our
arguments is using the properties of the measure established in the previous
section (in particular, Corollary~\ref{cor:highly-unbiasing-tail}) to get a
better bound on the probability that a disagreement percolates upwards. (In
\cite{MSW1} worst boundary conditions were assumed for the upward
coupling.)

We will establish a slightly different, but still equivalent to that stated in
Theorem~\ref{th:non-reconstruction}, form of the decay of correlations. For
every color $c$, let $\mu^{\downarrow}_{c,\ell}$ be the distribution of the
coloring of $L_\ell$ conditioned on the color at the root being $c$. For every
coloring $X$ of the leaves $L_\ell$ let $\mu^{\uparrow}_{X,\ell}$ be the
distribution of colors at the root conditioned on $X$ being the coloring of the
leaves. Finally, let $\mu^{\downarrow\uparrow}_{c,\ell}$ be the distribution of
colors at the root resulting from first choosing $X$ from
$\mu^{\downarrow}_{c,\ell}$ and then choosing a color for the root from
$\mu^{\uparrow}_{X,\ell}$.  Non-reconstruction is equivalent to the
dependence of $\mu^{\downarrow\uparrow}_{c,\ell}$ on~$c$ vanishing with~$\ell$.
We establish the decay of this dependence in the following theorem.
\begin{theorem}
\label{th:down-up-decay}
  Let $C>1$. There is a constant $\alpha(\eps)>0$ such that for all
  $\Delta > \Delta_0$, every $k > C\Delta/\ln\Delta$, and
  every pair of colors
  $c_1, c_2 \in  [k]$, for $\ell>\ell_0$,
  $$
  d_{TV} [\mu^{\downarrow\uparrow}_{c_1,\ell}, \mu^{\downarrow\uparrow}_{c_2,\ell}] \leq
  \Delta^{-\alpha\ell} .
  $$
\end{theorem}
Indeed, it is not too difficult to see that this theorem implies
Theorem~\ref{th:non-reconstruction}. (See Section~\ref{sec:appendix} for a proof of the
equivalence between the two forms of decay of correlations.)

We prove Theorem~\ref{th:down-up-decay} by coupling the two distributions.  The
coupling consists of two steps, the first is coupling the colorings of the
leaves (downward coupling) conditioned on the disagreement at the root, and the
second is coupling the color at the root (upward coupling) based on the
pair of colorings of the leaves achieved in the first step.

We start by bounding the average {\em Hamming} distance between the two coupled
colorings of the leaves resulting from the first step.  For two colorings $X,Y$,
let $H_\ell(X,Y)$ stand for the Hamming distance between $X(L_\ell)$ and
$Y(L_\ell)$, i.e., the number of leaves of the tree of depth~$\ell$ in which $X$
and $Y$ differ.
\begin{lemma}
\label{lem:downward-coupling}
  For every $\ell$ and any $c_1,c_2\in[k]$ there is a coupling
  $\nu=\nu^{\downarrow}_{c_1,c_2,\ell}$ of $\mu^{\downarrow}_{c_1,\ell}$ and
  $\mu^{\downarrow}_{c_2,\ell}$ such that
  $$
  \E_{(X,Y)\sim\nu} \left[H_\ell(X,Y)\right] =
  \left(\frac{\Delta}{k-1}\right)^{\ell}  \le \ln^\ell(\Delta).
  $$
\end{lemma}
\begin{proof}
  We let $\nu^{\downarrow}_{c_1,c_2,\ell}$ be the following recursive coupling:
  For $\ell \ge 1$, starting from the disagreement $(c_1,c_2)$ at the root,
  couple the spins at each of the children of the root independently according
  to the optimal coupling for each child. For each child vertex, if the resulting spins
  at this vertex agree then the whole subtree can be coupled to agree. If the
  spins disagree, say a $(c'_1, c'_2)$ disagreement, continue recursively using
  the coupling $\nu^{\downarrow}_{c'_1,c'_2,\ell-1}$, where for the base case
  the coupling $\nu^{\downarrow}_{c_1,c_2,0}$ is trivial. Notice that in our
  setting, the probability of a disagreement percolating to a given child is
  $\kappa = 1/(k-1)$.  To see this notice that the distribution at the child
  conditioned on the color of the root being~$c$ is uniform over
  $[k]\setminus\{c\}$. Thus, conditioned on a $(c_1,c_2)$-disagreement at the
  root, the two distributions at the child can be coupled such that the colors
  agree with probability $\frac{k-2}{k-1}$, i.e., that they agree whenever a
  color other than $c_1,c_2$ is chosen.  Clearly, the average Hamming distance
  at the leaves is $\kappa^{\ell}\Delta^{\ell} =
  \left(\frac{\Delta}{k-1}\right)^{\ell}$.
\end{proof}

The next ingredient we need is a bound on the distance between two conditional
distributions at the root in terms of the Hamming distance between the two
boundary colorings. We establish such a bound for boundary colorings that are
highly unbiasing. We say that a pair $(X,Y)$ of colorings of the leaves is {\em
  good} if both $(X,Y)$ are highly unbiasing.
\begin{lemma}
\label{lem:good-tv}
  For any given $\ell$ and any good pair $(X,Y)$ of colorings of the leaves
  $L_\ell$
  $$
  d_{TV} [\mu^{\uparrow}_{X,\ell} \,,\,
          \mu^{\uparrow}_{Y,\ell}] \;\le\;
  2H_\ell(X,Y) (2\Delta^{-\eps/2})^{(1-\eps)\ell}.
  $$
\end{lemma}
Let us complete the proof of Theorem~\ref{th:down-up-decay} assuming
Lemma~\ref{lem:good-tv}. For any pair $(X,Y)$ of colorings of the leaves
$L_\ell$ let $\nu^{\uparrow}_{X,Y,\ell}$ be the optimal coupling of
$\mu^{\uparrow}_{X,\ell}$ and $\mu^{\uparrow}_{Y,\ell}$, i.e., the probability
of disagreement between the pair of colors chosen in this coupling is exactly
$d_{TV} [\mu^{\uparrow}_{X,\ell} \,,\,
\mu^{\uparrow}_{Y,\ell} ]$.  We then set the coupling
$\nu^{\downarrow\uparrow}_{c_1,c_2,\ell}$ as follows.  First, choose a pair
$(X,Y)$ of colorings of the leaves from $\nu^{\downarrow}_{c_1,c_2,\ell}$; then,
choose a pair of colors at the root from the coupling
$\nu^{\uparrow}_{X,Y,\ell}$. Notice that
$\nu^{\downarrow\uparrow}_{c_1,c_2,\ell}$ is indeed a coupling of
$\mu^{\downarrow\uparrow}_{c_1,\ell}$ and $\mu^{\downarrow\uparrow}_{c_2,\ell}$.
We need to bound the probability that the two colors at the root chosen
according to this coupling disagree. Let $(c,c')$ be the pair of colors chosen from
this couplings. We have
\begin{eqnarray*}
\Prob[c \ne c'] &=& \sum_{X,Y} \nu^{\downarrow}_{c_1,c_2,\ell} (X,Y) \cdot
d_{TV}[\mu^{\uparrow}_{X,\ell} \,,\, \mu^{\uparrow}_{Y,\ell}]\\
&\le& \sum_{\text{good $(X,Y)$}} \nu^{\downarrow}_{c_1,c_2,\ell} (X,Y)
\cdot 2H_\ell(X,Y) (2\Delta^{-\eps/2})^{(1-\eps)\ell}
\\
& & \;+\;
\sum_{\text{bad $(X,Y)$}} \nu^{\downarrow}_{c_1,c_2,\ell} (X,Y)\cdot 1\\
&\le& 2(2\Delta^{-\eps/2})^{(1-\eps)\ell}
\sum_{X,Y} \nu^{\downarrow}_{c_1,c_2,\ell} (X,Y) \cdot H_\ell(X,Y)
\\
& & \;+\;
2\Prob_{X\sim\mu_\ell} [\text{$X$ is not highly unbiasing}]\\
&\le & 2(2\Delta^{-\eps/2})^{(1-\eps)\ell} \cdot \ln^\ell(\Delta) \;+\;
2e^{-\Delta^{\varepsilon \ell/3}}\\
&\le & \Delta^{-\eps\ell/4}
\end{eqnarray*}
assuming $\Delta \ge \Delta_0(\eps)$. Notice that in the first inequality above
we used Lemma~\ref{lem:good-tv}, in the second we used a simple union bound
together with the fact that the event of having a highly-unbiasing coloring of
the leaves is independent of the color at the root, and in the third inequality
we applied Lemma~\ref{lem:downward-coupling} to bound the first expression and
Corollary~\ref{cor:highly-unbiasing-tail} to bound the second. This completes
the proof of Theorem~\ref{th:down-up-decay} assuming Lemma~\ref{lem:good-tv}.

\subsection{Proof of Lemma~\ref{lem:good-tv}}
The proof of the lemma goes along similar lines to those in \cite{MSW1},
i.e., reducing to the case of a single disagreement using a sequence
of boundary
colorings where one leaf is changed at a time, and then using a coupling
argument to bound the probability that the disagreement at the leaf percolates
all the way up to the root. However, our argument differs in two points. The
first is that here we have to make sure that the intermediate boundary
colorings are
also valid, and more importantly, that they are highly-unbiasing. The
second difference is that we use the fact that the colorings are highly
unbiasing to get a better bound on the probability of a disagreement
percolating
upwards. (In our setting, if we do not assume highly unbiasing colorings at the
leaves then in the worst case the disagreement may percolate upwards with
probability~$1$.)

Let $(X,Y)$ be a pair of colorings of the leaves, $m=H(X,Y)$, and
$(d_1,\ldots,d_m)$ be an enumeration of the leaves in which $X$ and $Y$ differ.
We construct a sequence of colorings
$X=Z_0,Z_1,\ldots,Z_m,\ldots,Z_{2m}=Y$ going from $X$ to $Y$ by
un-conditioning the sites in which $X,Y$ differ one at a time and then adding the
condition according to $Y$ one site at a time. Formally, for $0\le i\le m$,
$Z_i(v)=\star$ if $v=d_j$ for some $1\le j \le i$, and otherwise $Z_i(v) =
X(v)$. For $i\ge m$, $Z_i(v)=\star$ if $v=d_j$ for some $1\le j \le 2m-i$, and
otherwise $Z_i(v)=Y$. It is easy to see that indeed for every $i\ge 0$,
$(Z_i,Z_{i+1})$ differ at exactly one leaf. We now show that each of the $Z_i$
is a valid highly-unbiasing coloring.
\begin{claim}
\label{prop:monotone}
  Let $X,X'$ be assignments from the leaves $L_\ell$ to $\{\star,1,\ldots, k\}$
  such that for all $v\in L_\ell$ we have $X'(v) \neq \star$ implies $X'(v) =
  X(v)$. If $X$ is a valid unbiasing coloring of $L_\ell$ then so is $X'$.
\end{claim}
\begin{proof}
  It is clear that if~$X$ is a valid coloring then so is~$X'$ because~$X'$ can
  be extended to $X$.  The unbiasing implication follows by induction on~$\ell$
  from the recursive definition of being unbiasing once we establish it for the
  base case of $\ell=1$. Indeed, for $\ell=1$, the number of unused colors in
  $X'$ is at least the number in $X$. Thus, if $X$ is unbiasing then so is $X'$.
\end{proof}
\begin{corollary}
\label{cor:partial-colorings-highly-unbiasing}
  Let $X,X'$ be assignments from the leaves $L_\ell$ to $\{\star,1,\ldots, k\}$
  such that for all $v\in L_\ell$ we have $X'(v) \neq \star$ implies $X'(v) =
  X(v)$. If $X$ is a valid highly-unbiasing coloring of $L_\ell$ then so is $X'$.
\end{corollary}

Now, by the triangle inequality, the proof of Lemma~\ref{lem:good-tv} will be
concluded once we show the following.
\begin{lemma}
\label{lem:single-disagreement}
  Let $Z,Z'$ be two highly-unbiasing colorings of $L_\ell$ that differ at
  exactly one vertex $v\in L_\ell$. Then,
  $$
  d_{TV} [\mu_\ell(\sigma(root) \,|\,\sigma(L)=Z) \,,\, \mu_\ell(\sigma(root)
  \,|\,\sigma(L)=Z')] \;\le\; (2\Delta^{-\eps/2})^{(1-\eps)\ell}.
  $$
\end{lemma}
The proof of the above lemma goes by an upward coupling argument similar to that
used in \cite{MSW1}, that is, we couple the spins along the path from~$v$ to the
root using a recursive coupling similar to that used for the downward coupling.
Specifically, given the disagreement at~$v$, we first couple the colors at the
immediate ancestor~$u$ (of~$v$) using an optimal coupling that minimizes the
probability of disagreement at~$u$. If the resulting spins at~$u$ agree, then we
couple the rest of the path with total agreement. If the spins at~$u$ disagree
then we continue recursively. The variation distance at the root is bounded by
the probability that the disagreement at the leaf~$v$ percolates along the path
from~$v$ all the way to the root. Thus, the important missing ingredient is a
bound on the probability of disagreement at a vertex~$u$ given a disagreement at
(exactly) one of its children (when the coloring of the leaves is~$Z$). This can
be done once we have a bound on the maximum probability of any color at~$u$, as
exemplified by the following proposition, which is taken from \cite{MSW2}.
\begin{proposition}
\label{prop:disagreement_prob}
  Let $\mu$ be a distribution of proper coloring of some graph according to some
  boundary condition, and let $\mu_u$ be the marginal of~$\mu$ at~$u$.  Consider
  now the distribution resulting from the same setting when adding a new
  neighbor to~$u$ (connected only to~$u$) and fixing the neighbor's color
  to~$c$. Let $\mu^c_u$ be the marginal of the latter distribution at~$u$. Then
  for any two colors $c_1,c_2$,
$$
d_{TV}[\mu^{c_2}_u, \mu^{c_1}_u] \;=\;
\max\left\{\mu^{c_2}_u(c_1)\,,\, \mu^{c_1}_u(c_2)\right\} \;=\;
\max\left\{\frac{\mu_u(c_1)}{1-\mu_u(c_2)}\,,\, \frac{\mu_u(c_2)}{1-\mu_u(c_1)}\right\}.
$$
In particular, if we let $p^{\max}_u = \max_c \mu_u(c)$ then
$d_{TV}[\mu^{c_2}_u, \mu^{c_1}_u] \le p^{\max}_u / (1-p^{\max}_u)$.
\end{proposition}
\begin{proof}
  We first notice that $\mu^c_u$ is $\mu_u$ conditioned on the resulting color
  not being~$c$, i.e., for every color $c'\ne c$, $\mu^c_u(c') =
  \mu_u(c')/(1-\mu_u(c))$. Assume without loss of generality that $\mu_u(c_1)\ge \mu_u(c_2)$. It
  follows that for every color $c'\ne c_1$, $\mu^{c_2}_u(c')\le
  \mu^{c_1}_u(c')$. Therefore, $d_{TV}[\mu^{c_2}_u, \mu^{c_1}_u] =
  \mu^{c_2}_u(c_1)$.
\end{proof}
Going back to our setting of the upward coupling,
Lemma~\ref{lem:single-disagreement} will follow once we show that for
highly-unbiasing $Z,Z'$, for every vertex~$w$ on the path from~$v$ to the root
which is at height $\ge \eps\ell$, given that the disagreement percolated
to~$w$, the probability that the disagreement percolates to the immediate
ancestor~$u$ of~$w$ is at most $2\Delta^{-\eps/2}$. Notice that once a color
at~$w$ is fixed, the distribution at~$u$ is independent of the rest of the
coloring of the subtree rooted at~$w$, and so that subtree (except for~$w$) can
be completely ignored. Thus, we are in position to use
Proposition~\ref{prop:disagreement_prob} to bound the probability of
disagreement at~$u$ given a disagreement at~$w$ by $p^{\max}_u /
(1-p^{\max}_u)$, where $p^{\max}_u$ in our setting is the maximum probability of
a color at~$u$ conditioned on~$Z$ at the leaves and when the edge $(u,w)$ is
removed from $T_\ell$ (i.e., the subtree rooted at~$w$ is removed from
$T_\ell$).

The fact that $Z$ is highly unbiasing
implies that the conditional probability $p^{\max}_u$ can be
bounded.

\begin{lemma}
\label{lem:unbiasing-to-bound-pmax}Let $u$ be a vertex of
  the tree $T_\ell$ at distance at least $\eps\ell$ from the leaves
  and $Z$ a highly unbiasing coloring of
  $L_\ell$. Let $p^{\max}_u$ be as defined above. Then,
  $$p^{\max}_u \leq \Delta^{-\eps/2}.$$
\end{lemma}

\begin{proof}
The proof is by the same arguments as in Lemma
  \ref{lem:unbiasing} so we don't repeat them, except to point out the
  following.
The probability considered in Lemma~\ref{lem:unbiasing} is
that of the subtree rooted at~$u$ conditioned on $Z_u$, i.e., this subtree is
disconnected from the rest of $T_\ell$. Here, the distribution considered is
when~$u$ is still connected to the rest of $T_\ell$, but with the subtree rooted
at~$w$ removed from the tree.  However, it is not too difficult to see that
if~$Z$ is highly unbiasing then the same bound on the maximum probability
applies to the distribution considered here.

Consider the orientation of $T_\ell$ (after removing the subtree
rooted at~$w$) so it is
rooted at~$u$. Since~$Z$ is a highly-unbiasing boundary
coloring of $T_\ell$ (in its original orientation) then at least $\Delta-1$ of
the children of~$u$ are unbiasing (those that are its children in the original
orientation), and hence~$Z$ is a highly unbiasing boundary condition for the
(irregular) newly-oriented tree rooted at~$u$. In the newly oriented tree, the
distribution considered in Lemma~\ref{lem:unbiasing} does correspond to that
considered here for $p^{\max}_u$.
\end{proof}

By the lemma, for every vertex $u$ at distance at least $\eps\ell$
from the leaves, $p^{\max}_u/(1-p^{\max}_u) \leq 2\Delta^{-\eps/2}$
and hence
$$
  d_{TV} [\mu_\ell(\sigma(root) \,|\,\sigma(L)=Z) \,,\, \mu_\ell(\sigma(root)
  \,|\,\sigma(L)=Z')] \;\le\; (2\Delta^{-\eps/2})^{(1-\eps)\ell},
  $$
completing the proof of
Lemma~\ref{lem:single-disagreement} and hence the proof of
Lemma~\ref{lem:good-tv}.

\subsection{Concentration}
\label{sec:concentration}
We now go on to prove Theorem~\ref{th:concentration}, i.e., establish that not
only does a typical boundary coloring yield a near uniform distribution of
colors at the root, but that the fraction of boundary colorings having a bias
at the root are extremely small. The main additional technical tool that we use here is a
reduction given in \cite{MSW1} from the tail bound stated in
Theorem~\ref{th:concentration} to an appropriate tail bound in the
coupling $\nu^{\downarrow}$.
\begin{lemma}
\label{lem:coupling-concentration}
  For any given $\ell$, $\delta>0$ and $A\ge 0$, if for every $c_1,c_2\in [k]$
  there exists a coupling~$\nu$ of $\mu^{\downarrow}_{c_1,\ell}$ and
  $\mu^{\downarrow}_{c_2,\ell}$ such that
  $$
  \Prob_{(X,Y)\sim\nu} \left[d_{TV}[\mu^{\uparrow}_{X,\ell} \,,\, \mu^{\uparrow}_{Y,\ell}]
  > \delta^3\right] \le A
  $$
  then for every $c\in [k]$,
  $$
  \Prob_{X\sim\mu_{\ell}} \left[|\mu^{\uparrow}_{X,\ell} (c) - 1/k| >
    2\delta\right] \le 2\left(e^{-1/\delta} + \frac{A}{\delta}\right).$$
\end{lemma}
\begin{proof}
  Arguments proving a similar statement to that of this lemma are given in
  Section~5.2 of \cite{MSW1}. For completeness, we repeat the proof in
  Section~\ref{sec:appendix}.
\end{proof}

As we show next, the coupling $\nu^{\downarrow}$ constructed above
satisfies the hypothesis of Lemma~\ref{lem:coupling-concentration} with the
necessary parameters. In fact, the only ingredient missing in establishing the
necessary tail bound for $\nu^{\downarrow}$ is a tail bound on the number of
disagreements in the coupling, which in turn follows from the (product)
percolation nature of the coupling.
\begin{lemma}
\label{lem:downward-coupling-tail}
  For every $\ell>\ell_0(\eps)$ and $c_1,c_2\in[k]$
  $$
  \Prob_{(X,Y)\sim\nu^{\downarrow}_{c_1,c_2,\ell}} \left[H_\ell(X,Y) >
    \Delta^{\eps\ell/8} \right] \leq e^{-\Delta^{\eps\ell/10}}.
  $$
\end{lemma}
\begin{proof}
  Recall from Lemma~\ref{lem:downward-coupling} that the average Hamming
  distance between $X,Y$ is $\le \ln^\ell\Delta$. The tail bound can then be
  established as done in, e.g., Lemma~5.5 of \cite{MSW1}.  For completeness, we
  give another proof in Section~\ref{sec:appendix}.
\end{proof}

Now, let us refine the notion of a pair of colorings $(X,Y)$ being good by
saying that a pair $(X,Y)$ of colorings of the leaves is {\em very good} if
$(X,Y)$ is good and $H_\ell(X,Y) \le \Delta^{\eps\ell/8}$. It follows from
Lemma~\ref{lem:good-tv} that for any very-good pair $(X,Y)$,
$$d_{TV}[\mu^{\uparrow}_{X,\ell}\,,\, \mu^{\uparrow}_{Y,\ell}] \;\le\;
\Delta^{-\eps \ell/8}.$$
On the other hand, by
Lemma~\ref{lem:downward-coupling-tail} and
Corollary~\ref{cor:highly-unbiasing-tail} (where the latter is used as in the
proof of Theorem~\ref{th:down-up-decay}) we have that
$$\Prob_{(X,Y)\sim\nu^{\downarrow}_{c_1,c_2,\ell}} \left[\text{$(X,Y)$
  is not very good}\right] \leq
  e^{-\Delta^{\eps\ell/10}} + 2e^{-\Delta^{\varepsilon \ell/3}} \le
  2e^{-\Delta^{\eps\ell/10}}.
$$
Thus, $\nu^{\downarrow}_{c_1,c_2,\ell}$ satisfies the hypothesis of
Lemma~\ref{lem:coupling-concentration} with $\delta=\Delta^{-3\ell}$ and
$A=2e^{-\Delta^{\eps\ell/10}}$, completing the proof of
Theorem~\ref{th:concentration}.

\section{Entropy constant for the block dynamics}
\label{sec:dynamics}
In this section we prove the rapid mixing result for the block dynamics as
stated in Theorem~\ref{th:dynamics}, i.e., we establish a lower bound on the
entropy constant associated with this dynamics when the block height is a large
enough constant. We rely on results established in \cite{MSW1}, where strong
concentration of the measure as established in Theorem~\ref{th:concentration}
above is shown to imply the necessary bound on the entropy constant of
the block
dynamics. We give here the high level motivation behind the kind of argument
made in \cite{MSW1} and some necessary adaptations to our setting,
while for the
main technical result we simply cite the relevant theorems from \cite{MSW1}.
Notice that while the arguments in \cite{MSW1} are written for the Ising model,
they hold for general models on trees. (We refer to \cite{DW-thesis} for the
general version of these arguments.)

To see the connection between strong concentration of the measure and the
entropy constant, recall that the dynamics makes heat-bath updates in a random
block, i.e., erases the local entropy in the updated block. Thus, the amount of
entropy erased in one step of the dynamics is the average (over blocks) of local
entropies, and in order to bound the entropy constant it is enough to show that
the global entropy of every function is not much more than the sum (over blocks)
of local entropies of the same function. The ability to express the global
entropy of every function as the sum of local entropies up to a constant factor
can be thought of as manifestation of locality in the equilibrium measure.
Indeed, the main technical result in \cite{MSW1} is that on trees, strong
concentration of the measure as in Theorem~\ref{th:concentration} is equivalent
to the ability to express entropy locally as above.

To make things formal, let $f$ be a function on the space of colorings $\Omega$.
The local entropy of $f$ in a subset $A$ of $T_n$ is $\mathbf{E}[\Ent (f
\,|\, \sigma_{T_n\setminus A})]$. Consider the block dynamics $\mathcal{M}_l$
and let $P_v$ be the transition matrix of the dynamics when the block to be
updated is the one rooted at $v$. Since the dynamics makes heat-bath updates of
the block $B_{v,\ell}$ then for any function $f$,
$P_v f (\sigma) = \E (f \,|\, \sigma_{T_n\setminus B_{v,\ell}})$, while $Pf = \frac{1}{N}
\sum_{v\in T_n} P_v f$. Now, for any subset $A$ of $T_n$ we have by
definition of entropy that
$$\Ent(f) = \Ent [ \E (f \,|\, \sigma_A)] + \E[ \Ent (f \,|\, \sigma_A)].$$
In particular, we have that
$$\Ent(f) - \Ent(P_v f) = \E[\Ent (f \,|\, \sigma_{T_n\setminus B_{v,\ell}})],$$
i.e., the effect of applying $P_v$ is to erase the local entropy of $f$ in $B_{v,\ell}$.
Let us denote the sum of local entropies as
$$\cE_\ell := \sum_{v\in T_n} \E[\Ent (f
\,|\, \sigma_{T_n\setminus B_{v,\ell}})].$$ Since entropy is
a convex functional then
$$\Ent(Pf) \;=\; \Ent \left(\frac{1}{N} \sum_{v\in T_n}
P_vf\right) \;\le\; \frac{1}{N} \sum_{v\in T_n} \Ent(P_v f).$$  In particular,
$$\Ent(f) - \Ent(Pf) \;\ge\; \frac{1}{N} \sum_{v\in T_n} [\Ent(f) -
\Ent(P_v f)] \;=\; \frac{1}{N}\cdot  \cE_\ell\enspace.$$
Thus, in order to prove the
lower bound on the entropy constant stated in Theorem~\ref{th:dynamics}, it is
enough to show that for an appropriate choice of the block height $\ell=\ell(\eps)$
(i.e., not depending on $n$),
\begin{equation}
\label{eqn:ent_decomposition}
\inf_f\frac{\cE_{\ell}(f)}{\Ent(f)} \ge \frac{1}{4}.
\end{equation}
Now, a key technical ingredient in~\cite{MSW1} is that strong concentration of the measure as established
in Theorem~\ref{th:concentration} implies~(\ref{eqn:ent_decomposition}). In order to state
the formal implication we need a couple of more definitions. For a vertex
$v\in T_n$, let $T_v$ be the subtree rooted at~$v$. For a color $c\in[k]$, let $\mu^c_{T_v}$
be the uniform distribution of colorings of $T_v$ conditioned on the parent of $v$ being
colored $c$. For $v$ whose distance from the leaves is at least~$\ell$, let $L_{v,\ell}$ be
the set of vertices at distance~$\ell$ below $v$. (Equivalently, $L_{v,\ell}$ is the bottom
boundary of $B_{v,\ell}$.)
The following theorem combines Theorems~3.4 and~5.3 of \cite{MSW1} (or Theorem~5.10 and
Lemma~5.18 in \cite{DW-thesis}, where the formulation is closer to the one here).
\begin{theorem}
\label{th:MSW}
For some absolute constant $\alpha^*$, the following implication is true for every choice of $\ell\ge 1$.
If for every vertex~$v$ whose distance from the leaves is at least~$\ell$ and for any two colors $c_1,c_2$,
$\mu\equiv\mu^{c_1}_{T_v}$ satisfies
\begin{equation}
\label{eqn:MSW}
\Prob_{\tau\sim\mu} \left[\left|\frac{\mu(\sigma(v)=c_2\,|\,\sigma(L_{v,\ell})=\tau(L_{v,\ell}))}{\mu(\sigma(v)=c_2)}-1\right| > \frac{1}{\alpha^* k^4(\ell+1)^2}\right] \le e^{-2\alpha^* k^4(\ell+1)^2}
\end{equation}
then for every function $f\ge 0$, $\Ent(f) \le 4 \cE_{\ell}$.
\end{theorem}
We conclude by noticing that from Theorem~\ref{th:concentration} (see the
following Remark \ref{rem:parent-color} for a proviso) we get
that~(\ref{eqn:MSW})  (and therefore~(\ref{eqn:ent_decomposition}))
holds for an appropriately chosen $\ell$, e.g., $\ell \geq \frac{10\alpha^*}{\min\{\alpha,\alpha'\}}$ suffices where $\alpha^*$ is the absolute constant from Theorem
\ref{th:MSW}, and $\alpha=\alpha(\eps)$ and $\alpha'=\alpha'(\eps)$
are from Theorem \ref{th:concentration}.

\begin{remark} 
\label{rem:parent-color}
  Notice that Theorem~\ref{th:MSW} requires us to prove the strong concentration of the
  measure stated in Theorem~\ref{th:concentration} in a slightly more general scenario.
  Specifically, we need to show similar bounds for a random boundary coloring of
  $B_{v,\ell}$ conditioned on an arbitrary color at the parent of~$v$ and when
  this boundary coloring is chosen from the measure on $T_n$. Notice, however,
  that the distribution of colorings at the boundary of $B_{v,\ell}$ under the
  distribution of colorings of~$T_n$ is exactly the same as the distribution of
  colorings of the leaves of the smaller tree in which the vertices on the
  boundary of $B_{v,\ell}$ are the leaves. As for fixing an arbitrary color
  above $v$, it is straightforward that the concentration stated in
  Theorem~\ref{th:concentration} holds in this scenario as well. (In this case
  the distribution at the root is uniform over $k-1$ colors so the statement of
  the theorem is modified with $\frac{1}{k-1}$ replacing $\frac{1}{k}$.)
\end{remark}
\begin{remark}
  In \cite{MSW1} (and \cite{MSW2,DW-thesis}), the bound
  in~(\ref{eqn:ent_decomposition}) is used to bound the mixing time of the
  single-site dynamics. This is done by introducing an extra factor that
  expresses the decomposition of the entropy in a block (under an arbitrary
  boundary condition) into the local entropy (or variance) coming from single
  sites in the block. This factor is bounded by the mixing time of the
  single-site dynamics in the block, which in the setting of
  \cite{MSW1,MSW2,DW-thesis} is a constant depending on the block size.
  However, in our setting there are valid boundary colorings of the block for
  which the single-site dynamics is disconnected (i.e., infinite mixing time) so
  we cannot apply the same reasoning and suffice with a bound on the mixing time
  of the block dynamics.
\end{remark}

\section{Leftover Proofs}
\label{sec:appendix}
Here we provide the remaining proofs.
Most of the arguments given
in this section already appear in the references mentioned in the main text, but we
give them here for completeness and in order for the proofs to correspond
exactly to the formulations of the statements made in this paper.

We start by showing that the form of decay of correlations stated in
Theorem~\ref{th:down-up-decay} is indeed equivalent to that stated in the main
Theorem~\ref{th:non-reconstruction}. Let $\alpha_{c,\ell} = \mathbf{E}_{X\sim
\mu_\ell} [|\mu^\uparrow_{X,\ell}(c) - 1/k|]$ and $\beta_{c,\ell} =
|\mu^{\downarrow\uparrow}_{c,\ell} (c) - 1/k|$. Notice that
Theorem~\ref{th:down-up-decay} establishes that for every $c$ and $\ell$
$\beta_{c,\ell} = \Delta^{-\Omega(\ell)}$, while for
Theorem~\ref{th:non-reconstruction} we need to show that for every $c$ and
$\ell$, $\alpha_{c,\ell} = \Delta^{-\Omega(\ell)}$. The following proposition shows
that the two statements are indeed equivalent.
\begin{proposition}
\label{prop:equiv_decay}
For every $c$ and $\ell$, $\frac{\beta_{c,\ell}}{k-1} \le \alpha_{c,\ell} \le \sqrt{\beta_{c,\ell}}$.
\end{proposition}
\begin{proof}
  Let $g_{c,\ell}(\cdot)$ be the function on coloring of the leaves $L_\ell$
  defined as follows:
  $$
  g_{c,\ell} (X) = \frac {\mu^{\downarrow}_{c,\ell}(X)}{\mu_\ell(X)} = \frac
  {\mu^{\uparrow}_{X,\ell}(c)}{\mu_\ell(\sigma_{root} = c)} = k\cdot
  \mu^{\uparrow}_{X,\ell}(c).
  $$
  Notice that $\E_{X\sim\mu_\ell} [g_c(X)] = 1$. Furthermore, $\alpha_{c,\ell} =
  \frac{1}{k}\E_{X\sim\mu_\ell} [|g_c(X) - 1|]$, and $\beta_{c,\ell} =
  \frac{1}{k}\E_{X\sim\mu_\ell} [|g_c(X)|^2 - 1] = \frac{1}{k}\E_{X\sim\mu_\ell}
  [|g_c(X) - 1|^2]$. The first inequality of the proposition then follows from
  the fact that $\|g_c\|_{\infty} \le k$, while the second inequality follows
  from Cauchy-Schwartz.
\end{proof}

We now provide the omitted proofs of the two lemmas in
Section~\ref{sec:concentration}.
\begin{proof}[Proof of Lemma~\ref{lem:coupling-concentration}]
  First, notice that since $\mu_\ell(\sigma(L_\ell))$ is a convex combination of
  $\mu^{\downarrow}_{c,\ell}$ we can assume w.l.o.g. that for every color~$c$ we
  have a coupling of $\mu^{\downarrow}_{c,\ell}$ and $\mu_\ell$ with
  properties as in the hypothesis of the lemma.  Using the same definition of
  $g_{c,\ell}$ as in Proposition~\ref{prop:equiv_decay}, for any $\alpha\ge 0$ we
  then have that
  $$
  \Prob_{X\sim\mu^{\downarrow}_{c,\ell}} [g_{c,\ell}(X) \ge 1 + \alpha] \;\le\;
  \Prob_{X\sim\mu_\ell} [g_{c,\ell}(X) \ge 1 + \alpha - k\delta^3] + A.
  $$
  On the other hand, by the definition of $g_{c,\ell}$,
  $$
  \Prob_{X\sim\mu^{\downarrow}_{c,\ell}} [g_{c,\ell}(X) \ge 1 + \alpha] \;\ge\;
  (1+\alpha)\Prob_{X\sim\mu_\ell} [g_{c,\ell}(X) \ge 1 + \alpha].
  $$
  Combining the above two inequalities we get that
  $$
  \Prob_{X\sim\mu_\ell} [g_{c,\ell}(X) \ge 1 + \alpha] \;\le\;
  \left(\frac{1}{1+\alpha}\right)\left(\Prob_{X\sim\mu_\ell} [g_{c,\ell}(X)
    \ge 1 + \alpha - k\delta^3] + A\right).
  $$
  In particular, for every non-negative integer~$m$, if we apply the above
  inequality $m+1$ times, increasing $\alpha$ by $k\delta^3$ each time, we get
  that
  $$
  \Prob_{X\sim\mu_\ell} [g_{c,\ell}(X) \ge 1 + \alpha + mk\delta^3] \;\le\;
  \left(1+\alpha\right)^{-(m+1)} \;+\; \frac{A}{k\alpha} \;\le\;
  e^{-\alpha(m+1)/(1+\alpha)} + \frac{A}{\alpha} \enspace .
  $$
  Now, by setting $\alpha=k\delta$ and $m=\lfloor(1/\delta)^2\rfloor$ we get that
  $$
  \Prob_{X\sim\mu_\ell} [g_{c,\ell}(X) \ge 1 + 2k\delta] \;\le\;
  e^{-\delta(m+1)} + \frac{A}{\delta} \enspace .
  $$
  Since the event $g_{c,\ell}(X) \ge 1 + 2k\delta$ is equivalent to
  $\mu^{\uparrow}_{X,\ell}(c) - 1/k \ge 2\delta$, this establishes the bound on
  the positive side of the tail in Lemma~\ref{lem:coupling-concentration}. The
  negative side of the tail is established along similar lines once we notice
  that from the hypothesis we have for any $\alpha\ge 0$
  $$
  \Prob_{X\sim\mu_\ell} [g_{c,\ell}(X) \le 1 - \alpha - k\delta^3] \;\le\;
  \Prob_{X\sim\mu^{\downarrow}_{c,\ell}} [g_{c,\ell}(X) \le 1 - \alpha] + A
  $$
  and by definition of $g_{c,\ell}$,
  $$
  \Prob_{X\sim\mu^{\downarrow}_{c,\ell}} [g_{c,\ell}(X) \le 1-\alpha] \;\le\;
  (1-\alpha)\Prob_{X\sim\mu_\ell} [g_{c,\ell}(X) \le 1-\alpha].
  $$
\end{proof}
\begin{proof}[Proof of Lemma~\ref{lem:downward-coupling-tail}]
  Recall that the coupling $\nu^{\downarrow}$ is constructed recursively such
  that if the spin at a vertex~$v$ agree then all of its children are coupled
  with agreement with probability $1$, while if there is a disagreement at~$v$
  then its children are coupled independently where the probability of
  disagreement in each of the $\Delta$ children is $1/(k-1)$.  Thus, if we let
  $D_i$ be the random variable counting the number of disagreements at level~$i$
  in this coupling then $D_0 = 1$ and $D_{i+1}$ is distributed as
  $\mathrm{Bin}(\Delta D_i, \frac{1}{k-1})$. We need to establish an upper bound
  on $\Prob(D_\ell \ge \Delta^{\eps\ell/8})$. Let $i_0 = \eps\ell/9$.  For every
  $i_0 \le i \le \ell$ we say that a failure occurred at level~$i$ if $D_i >
  (\frac{3\ln\Delta}{2})^{i-i_0} \cdot \Delta^{i_0}$. Now, notice that given
  $D_i$, $\E[D_{i+1}] = \frac{\Delta}{k-1} D_i \le D_i \ln \Delta$. Thus, by a
  Chernoff bound,
  for every $i\ge i_0$,
  $$
  \Prob[\text{ failure in level $i+1$} \,|\, \text{not failure in level $i$}]
  \;\le\; \exp \left(-\frac{\ln\Delta}{16} \cdot
    \left(\frac{3\ln\Delta}{2}\right)^{i-i_0} \Delta^{i_0} \right).
  $$
  Now, since $\Prob[\text{ failure in level $i+1$}] \le   \Prob[\text{ failure in
  level $i+1$} \,|\, \text{not failure in level $i$}] + \Prob[\text{ failure in
  level $i$}]$ and since $\Prob[\text{ failure in level $i_0$}] = 0$, we have by
  induction that
\begin{eqnarray*}
\Prob[\text{ failure in level $\ell$}] &\le& \sum_{i=i_0+1}^{\ell}
  \Prob[\text{ failure in
  level $i$} \,|\, \text{not failure in level $i-1$}]\\
&\le&
  2\exp\left(\frac{-\Delta^{i_0}}{16}\right)\\
&\le& \exp\left(-\Delta^{\eps\ell/10}\right).
\end{eqnarray*}
Notice that the event of failure in level $\ell$ includes the event
that $D_\ell \ge
\Delta^{\eps \ell /8}$ so we are done.
\end{proof}

\section*{Acknowledgments} The authors would like to thank Fabio
  Martinelli, Elchanan Mossel, Alistair Sinclair, Allan Sly and Prasad
  Tetali for helpful discussions on the topic.

\end{document}